\documentclass[12pt,reqno]{amsart}
\usepackage{amssymb,amsthm} 
\usepackage[utf8]{inputenc}

\advance\hoffset by -0.5truein
\let\<\langle
\let\>\rangle
\def\oo#1{\mathrel { {}_{(#1)}}}
\def\Var{\mathop {\fam0 Var}\nolimits}
\def\Ker{\mathop {\fam0 Ker}\nolimits}
\def\Cur{\mathop {\fam0 Cur}\nolimits}
\def\Chom{\mathop {\fam0 Chom}\nolimits}
\def\Cend{\mathop {\fam0 Cend}\nolimits}
\def\gc{\mathop {\fam0 gc}\nolimits}
\def\Pois{\mathrm{Pois}}

\def\GD {\mathop {\fam0 GD }\nolimits}

\def\Com{\mathrm{Com}}
\def\Nov{\mathrm{Nov}}

\def\Der{\mathrm{Der}}
\def\wt{\mathop {\fam 0 wt}\nolimits}

\newtheorem{theorem}{Theorem}
\newtheorem{lemma}{Lemma}
\newtheorem{proposition}{Proposition}
\newtheorem{corollary}{Corollary}

\theoremstyle{definition}
\newtheorem{definition}{Definition}

\newtheorem{example}{Example}

\title[Quadratic Lie conformal superalgebras]{Quadratic Lie conformal superalgebras related to Novikov superalgebras}

\thanks{The work is supported by Mathematical Center in Akademgorodok under agreement 
No.~075-15-2019-1613 with the Ministry of Science and Higher Education of the Russian Federation}

\author[P.S. Kolesnikov, R.A. Kozlov, A.S. Panasenko]{P.S. Kolesnikov$^{1)}$, R.A. Kozlov$^{1)2)}$, A.S. Panasenko$^{1)2)}$}

\address{${}^{1)}$ Sobolev Institute of Mathematics, Novosibirsk, Russia}
\address{${}^{2)}$ Novosibirsk State University, Novosibirsk, Russia}

\begin{document}
\begin{abstract}
We study quadratic Lie conformal superalgebras associated with No\-vikov superalgebras. 
For every Novikov superalgebra $(V,\circ)$, we construct an enveloping differential Poisson superalgebra $U(V)$ 
with a derivation $d$
such that $u\circ v = ud(v)$ and $\{u,v\} = u\circ v - (-1)^{|u||v|} v\circ u$ for $u,v\in V$.
The latter means that the commutator Gelfand--Dorfman superalgebra of $V$ is special.
Next, we prove that every quadratic Lie conformal superalgebra constructed on a finite-dimensional special 
Gelfand--Dorfman superalgebra has a finite faithful conformal representation. 
This statement is a step toward a solution of the following open problem: 
whether a finite Lie conformal (super)algebra has a finite faithful conformal representation.
\end{abstract}


\maketitle

\section{Introduction}

Novikov algebras appeared in \cite{GD83} as a class of algebras giving rise to 
Hamiltonian operators in the formal calculus of variations. Independently, 
these algebras were introduced in \cite{BN85} as a tool for studying linear Poisson brackets of hydrodynamic type.
The study of structure theory of Novikov algebras was initiated in \cite{Zelm},  
significant progress in this direction was obtained in \cite{BaiDeng01, BdG2013, Xu1996, Xu2001}.

A class of more complicated structures called Gelfand--Dorfman bialgebras \cite{Xu2000} was also introduced in \cite{GD83}
as a source of Hamiltonian operators. A Gelfand--Dorfman bialgebra is a linear space with 
two bilinear operations $(\cdot \circ \cdot)$ and $[\cdot, \cdot]$, where 
$(\cdot \circ \cdot)$ is a Novikov product (left symmetric, right commutative), i.e., 
\begin{gather}
(x_1\circ x_2)\circ x_3-x_1\circ (x_2\circ x_3)=(x_2\circ x_1)\circ x_3 - x_2\circ (x_1\circ x_3), \label{eq:LSymm} \\
(x_1\circ x_2)\circ x_3 = (x_1\circ x_3)\circ x_2, \label{eq:RComm}
\end{gather}
$[\cdot, \cdot]$ is a Lie product, and the following compatibility relation holds:
\begin{equation}\label{eq:GD-1} 
[x_1,x_2\circ x_3]-[x_3,x_2\circ x_1]+[x_2,x_1]\circ x_3 - [x_2,x_3]\circ x_1 - x_2\circ[x_1,x_3]=0.
\end{equation}
In order to avoid a confusion with the well-known notion of a bialgebra as an algebra equipped
with a coproduct, we will use the term {\em GD-algebra} (\cite{KSO2019}) 
for a Gelfand---Dorfman bialgebra.

Novikov algebras and GD-algebras play an important role in the combinatorics 
of differential algebras (see \cite[Theorem~7]{KSO2019}). 
Namely, 
the identities that hold for the operations 
\[
 a\succ b  = d(a)\cdot b,\quad a\prec b = a\cdot d(b),
\]
$a,b\in A$, with
$A$ being an algebra in a multilinear variety $\Var $,
may be calculated by means of the Manin white 
product \cite{GK94} of the operad $\Var $ and 
the operads $\Nov $ and $\GD^!$, the Koszul dual of~$\GD$.

Every differential Poisson algebra gives rise to a GD-algebra \cite{Xu2002}: 
if $P$ is a commutative differential algebra with a derivation $d$ 
equipped with a Poisson bracket $\{\cdot,\cdot\}$ such that $d$ acts as a 
derivation relative to $\{\cdot, \cdot\}$ then 
$P$ is a GD-algebra with operations $x\circ y =zd(y)$, $[x,y]=\{x,y\}$, $x,y\in P$. 
GD-algebras that embed into Poisson differential algebras in this way are called 
special \cite{KSO2019}. Not all GD-algebras are special, 
a series of necessary conditions for a GD-algebra to be special 
was found in the last paper. Below, we will state some precise examples of non-special GD-algebras. 

It was noted in \cite[Remark~6.3]{GD83} that if we enrich a Novikov algebra  $(V,\circ)$ with 
the operation $[x,y] = x\circ y-y\circ x$, $x,y\in V$, 
then we obtain a GD-algebra. It is not hard to check that such a system obtained from 
a Novikov algebra meets necessary conditions of speciality found in \cite{KSO2019}.
In Section \ref{sec:NovikovEnv}, we prove that all GD-superalgebras arising from  Novikov superalgebras relative 
to (super-)commutator 
are indeed special, and we construct an enveloping differential Poisson superalgebra for every Novikov algebra.

It turned out (see \cite{HongWu17, Xu2000}) that GD-algebras are in one-to-one correspondence with 
quadratic Lie conformal algebras. The latter structures appeared in \cite{Kac1998}
as a tool in the study of vertex operator algebras. 
Conformal algebras and their generalizations (pseudo-algebras) also turn to be useful 
for the classification of Poisson brackets of hydrodynamic type \cite{BDK}. 
A conformal algebra is a module $C$ over the polynomial algebra $H=\Bbbk [\partial ]$ 
equipped with a ``multi-valued'' operation $C\otimes C\to C[\lambda ]$, 
i.e., a product of two elements from $C$ is a polynomial in a formal variable 
$\lambda $ with coefficients in~$C$. The axioms defining a conformal algebra 
are stated in Section~\ref{sec:Conformal}.

One of the most intriguing questions in the theory of conformal algebras 
is motivated by the Ado theorem: Does a Lie conformal algebra 
which is a finitely generated free module over $H$ have a faithful 
representation on a finitely generated free $H$-module?
(The condition of freeness is necessary: in conformal algebras and their modules, 
every torsion element belongs to the corresponding annihilator.)
In the case of positive answer, we may faithfully represent every 
Lie conformal algebra with polynomial matrices in $\Cend_n$, see \cite{BKL2003}.
By now, the most general result in this direction says that 
if a Lie conformal algebra as above has a Levi decomposition (i.e., 
its solvable radical splits) then it has a finite faithful representation (FFR).

Therefore, one may observe close relations between Novikov algebras, Poisson algebras, and conformal algebras.
In Section~\ref{sec:Poisson}, we prove that every quadratic Lie conformal superalgebra 
obtained from a special GD-algebra has a FFR. As an application, every quadratic Lie conformal superalgebra 
obtained from a Novikov superalgebra has a FFR.

\section{Gelfand--Dorfman superalgebras}\label{sec:NovikovEnv}

A $\mathbb Z_2$-graded space $V=V_0\oplus V_1$ with two bilinear operations 
$(\cdot\circ\cdot ), [\cdot,\cdot ]: V\otimes V\to V$ which respect the grading 
is said to be a  GD-superalgebra if 
\begin{itemize}
 \item $V$ is a Novikov superalgebra relative to $ (\cdot\circ \cdot)$; 
 \item $V$ is a Lie superalgebra relative to $ [\cdot,\cdot]$; 
 \item for every homogeneous $a,b,c\in V$
 \begin{equation}\label{eq:GD-1super} 
[a\circ b,c] - a\circ[b,c] + [a,b]\circ c + (-1)^{|a||b|}[b,a\circ c] - (-1)^{|b||c|} [a,c]\circ b = 0 .
\end{equation}
\end{itemize}

A series of examples of GD-superalgebras may be obtained from differential Poisson superalgebras. 
Let $P=P_0\oplus P_1$ be an associative and commutative superalgebra equipped with an operation $\{\cdot ,\cdot \}$
such that 
\begin{gather}
 \{P_i,P_j\}\subseteq P_{(i+j)\pmod 2}, \quad \{a,b\}=-(-1)^{|a||b|}\{b,a\}, \nonumber \\
\{a,\{b,c\}\} - (-1)^{|a||b|} \{b,\{a,c\}\} = \{\{a,b\},c\},  \label{eq:sJacobi}\\
\{a,bc\} = \{a,b\}c +(-1)^{|a||b|} b\{a,c\}, \label{eq:sLeibniz}
\end{gather} 
for all homogeneous $a,b,c\in P$.
Then $P$ is called a Poisson superalgebra. 

Suppose a Poisson superalgebra $P$ has an (even) derivation $d:P\to P$, i.e., 
\[
 d(P_i)\subseteq P_i,\quad d(ab)=d(a)b+ad(b),\quad d\{a,b\} = \{d(a),b\} + \{a,d(b)\},
\]
for all $a,b\in P$. 
Then the same space $P$ equipped with 
\[
 a\circ b = ad(b),\quad [a,b]= \{a,b\}
\]
is a GD-superalgebra \cite[Theorem~3.2]{Xu2000} denoted $P^{(d)}$.

For a GD-superalgebra $V$, we say $V$ is {\em special} if there exists a 
Poisson superalgebra $P$ with a derivation $d$ such that $V\subseteq P^{(d)}$.
Non-special GD-superalgebras are said to be {\em exceptional}.

Exceptional GD-superalgebras exist: it was shown in \cite{KSO2019} by 
means of implicit computational arguments. Let us state below
an explicit example of a 3-dimensional exceptional GD-algebra.

\begin{example}\label{exmp:Heisenberg3}
Let $V$ be a 3-dimensional space with a basis $\{x,y,z\}$ 
equipped with a Lie algebra structure 
\[
 [x,y]=z,\quad [x,z]=[y,z]=0
\]
(Heisenberg Lie algebra). It is straightforward to check that the operation $(\cdot\circ\cdot)$
on $V$ given by 
\[
 \begin{gathered}
 x\circ x = x-y,\quad y\circ x =  - x\circ y = y ,  \\
 x\circ z=z\circ x = y\circ z = z\circ y = y\circ y = z\circ z = 0
 \end{gathered}
\]
turns $V$ into a GD-algebra. If we assume $V$ to be special, then there exists a 
Poisson differential algebra $P$ such that $V\subset P^{(d)}$. Consider 
the expression $\{x,x'\}xx' \in P$ (here $x'$ stands for $d(x)$).
On the one hand,
\begin{multline}\nonumber
 \{x,x'\}(xx') = \{x,x'\}(x\circ x) = \{x,(x\circ x)x'\} -  \{x,x\circ x\}x' \\
 = [x,(x\circ x)\circ x] - [x,x\circ x]\circ x = -2z;
\end{multline}
on the other hand,
\[
\{x,x'\}(xx')= (\{x,x'\}x)x' = \{x,xx'\}x' 
= [x,x\circ x]\circ x = 0.
\]
Hence, $V$ cannot be embedded into a differential Poisson algebra.
\end{example}

Another series of examples was proposed in \cite{GD83} for the non-graded case. 
In general, if $V$ is a Novikov superalgebra then the operation 
\begin{equation}
 [a,b] = (a\circ b) - (-1)^{|a||b|} (b\circ a) , \quad a,b\in V_0\cup V_1,
\end{equation}
turns $V$ into a GD-superalgebra denoted $V^{(-)}$.
It was conjectured in \cite{KSO2019} that 
for every Novikov algebra the commutator GD-algebra $V^{(-)}$ is special.
In this section, we prove the conjecture in the $\mathbb Z_2$-graded setting.

In the following, we will need the embedding theorem proved in \cite{BCZ-2017}
for Novikov algebras and then in \cite{ZCB2019} for Novikov superalgebras
(in these two papers, the term ``Gelfand--Dorfman--Novikov algebra'' is used 
for what we call ``Novikov algebra'' following the common terminology proposed in \cite{Osborn}).

\begin{theorem}[{\cite[Theorem~3]{BCZ-2017}, \cite[Theorem 3.8]{ZCB2019}}]\label{thm:NovToComDer}
 For every Novikov (super)algebra $(V,\circ)$ 
 there exists an associative and (super)commutative 
 algebra $A$ with an even derivation $d$
 such that $V\subset A$ and 
 $u\circ v = ud(v)$ for $u,v\in V$.
\end{theorem}

In particular, the universal enveloping differential algebra of a given 
Novikov superalgebra $V=V_0\oplus V_1$ may be constructed as follows.
Let $X=X_0\cup X_1$ be a linear basis of $V$, where $X_i$ is a basis of $V_i$, $i=0,1$. 
Denote by $s\Com\Der\<X_0\cup X_1, d\>$ the free associative supercommutative differential algebra
generated by even variables $X_0$ and odd variables $X_1$ with an even derivation~$d$. 
Apparently, $s\Com\Der\<X_0\cup X_1, d\>\simeq \Bbbk [d^\omega X_0]\otimes \bigwedge (\Bbbk d^\omega X_1)$, 
where 
\[
 d^\omega X_i  = \{ d^n(x) \mid n\ge 0, \ x\in X_i \}.
\]
Consider the (differential) ideal $I_V$ of $s\Com\Der\<X_0\cup X_1, d\>$ generated by 
$u(x,y) = xd(y)-x\circ y$, $x,y\in X$ (here $x\circ y$ stands for the linear form in $\Bbbk X$ representing 
the Novikov product in $V$). As a non-differential ideal of 
$\Bbbk [d^\omega X_0]\otimes \bigwedge (\Bbbk d^\omega X_1)$,
$I_V$ is generated by all derivatives of $u(x,y)$. Theorem \ref{thm:NovToComDer} implies $I_V\cap \Bbbk X = 0$.

In order to define a (super) Poisson bracket $\{\cdot,\cdot\}$ on $s\Com\Der\<X_0\cup X_1, d\>$
it is enough to determine polynomials $\{d^n(x),d^m(y)\}$ for $x,y\in X$, $n,m\ge 0$, 
and then extend the bracket in a unique way to the entire $s\Com\Der\<X_0\cup X_1, d\>$ by the Leibniz rule
\eqref{eq:sLeibniz}. 
If the bracket respects $\mathbb Z_2$-grading and is (super) anti-commutative on the generators from $d^{\omega }X$ 
then so is 
its extension. 
In order to simplify notations, let us denote $d^n(x)$ by $x^{(n)}$ for $x\in X$, $n\ge 0$.

\begin{lemma}\label{lem:BrJacobi}
Let $\{\cdot, \cdot \}$ be a bracket on $s\Com\Der\<X_0\cup X_1, d\>$ obtained by 
expanding 
\begin{equation}\label{eq:NovPoisson_s_Bracket}
\{x^{(m)},y^{(n)}\}=(n-1)x^{(m+1)}y^{(n)}-(m-1)x^{(m)}y^{(n+1)}, \quad x,y\in X,\ n,m\ge 0.
\end{equation}
Then $\{\cdot,\cdot\}$ satisfies the Jacobi identity \eqref{eq:sJacobi}.
\end{lemma}

\begin{proof}
For $x,y,z\in X$, $n,m,k\ge 0$, evaluate 
\begin{multline}\label{eq:3-bracket}
\{x^{(m)}\{y^{(n)},z^{(k)}\}\}
=\{x^{(m)},(k-1)y^{(n+1)}z^{(k)}-(n-1)y^{(n)}z^{(k+1)}\} \\
=(k-1)(\{x^{(m)},y^{(n+1)}z^{(k)}+(-1)^{|x||y|}y^{(n+1)}\{x^{(m)},z^{(k)}\}\}) \\
-(n-1)(\{x^{(m)},y^{(n)}\}z^{(k+1)}+(-1)^{|x||y|}y^{(n)}\{x^{(m)},z^{(k+1)}\}) \\
=(k-1)\big((nx^{(m+1)}y^{(n+1)}-(m-1)x^{(m)}y^{(n+2)})z^{(k)} \\
+y^{(n+1)}(-1)^{|x||y|}((k-1)x^{(m+1)}z^{(k)}-(m-1)x^{(m)}z^{(k+1)})\big) \\
-(n-1)\big(((n-1)x^{(m+1)}y^{(n)}-(m-1)x^{(m)}y^{(n+1)})z^{(k)} \\
-(-1)^{|x||y|}y^{(n)}(kx^{(m)}z^{(k+1)}-(m-1)x^{(m)}z^{(k+2)})\big)\\
=(k-1)(n+k-1)x^{(m+1)}y^{(n+1)}z^{(k)}-(k-1)(m-1)x^{(m)}y^{(n+2)}z^{(k)} \\
+(m-1)(n-k)x^{(m)}y^{(n+1)}z^{(k)}-(n-1)(k+n-1)x^{(m+1)}y^{(n)}z^{(k+1)} \\
+(n-1)(m-1)x^{(m)}y^{(n)}z^{(k+2)}.
\end{multline}
The Jacobi identity \eqref{eq:sJacobi} is equivalent to 
\begin{multline}\nonumber
(-1)^{|x||z|}\{x^{(m)},\{y^{(n)},z^{(k)}\}\}+(-1)^{|y||x|}\{y^{(n)},\{z^{(k)},x^{(m)}\}\} \\
+(-1)^{|y||z|}\{z^{(k)},\{x^{(m)},y^{(n)}\}\}=0
\end{multline}
which is easy to check by the cyclic permutation of variables in \eqref{eq:3-bracket}.
\end{proof}

\begin{lemma}\label{lem:Der_bracket}
The operation $d$ on $s\Com\Der\<X_0\cup X_1, d\>$
is a derivation relative to the bracket from Lemma~\ref{lem:BrJacobi}.
\end{lemma}

\begin{proof}
It is enough to check that 
\begin{multline*}
d\{x^{(m)}, y^{(n)}\}=(n-1)(x^{(m+2)}y^{(n)}+x^{(m+1)}y^{(n+1)})
-(m-1)(x^{(m)}x^{(n+2)}+x^{(m+1)}y^{(n+1)}) \\
=(n-1)x^{(m+2)}y^{(n)}-mx^{(m+1)}y^{(n+1)}+x^{(m+1)}y^{(n+1)} +nx^{(m+1)}y^{(n+1)}\\
-(m-1)x^{(m+1)}y^{(n+1)}-x^{(m+1)}y^{(n+1)} 
  =\{x^{(m+1)},y^{(n)}\}+\{x^{(m)},y^{(n+1)}\} \\
= \{dx^{(m)},y^{(n)}\} + \{x^{(m)},dy^{(n)}\}
\end{multline*}
for $x,y\in X$, $n,m\ge 0$. 
Since the bracket on the entire $s\Com\Der\<X_0\cup X_1, d\>$
is calculated via the Leibniz rule, we have 
\[
d\{f,g\} = \{d(f),g\} + \{f,d(g)\}
\]
for all $f,g\in s\Com\Der\<X_0\cup X_1, d\>$.
\end{proof}

\begin{lemma}\label{lem:bracket_Ideal}
The ideal $I_V$ 
is invariant under all operations $\{f,\cdot\}$, $f\in s\Com\Der\<X_0\cup X_1, d\>$,
where $\{\cdot,\cdot \}$ is the bracket from Lemma~\ref{lem:BrJacobi}.
\end{lemma}

\begin{proof}
Again, since $\{f,\cdot\}$ is defined via the Leibniz rule, it is enough to consider $f=z^{(n)}$, $z\in X$, $n\ge 0$.
Moreover, Lemma~\ref{lem:Der_bracket} implies
\[
\{z^{(n)}, \cdot \} = d\{z^{(n-1)}, \cdot\} - \{z^{(n-1)}, d(\cdot ) \}, 
\]
so it is enough to check the invariance of $I_V$ under $\{z,\cdot\}$, where $z=z^{(0)}$.
The Leibniz rule \eqref{eq:sLeibniz} and Lemma~\ref{lem:Der_bracket} show that 
it is enough to verify $\{z,u(x,y)\} \in I_V$ for all $x,y,z\in X$. Indeed, 
\begin{multline*} 
\{z, u(x,y) \}=\{z, xy'\}-\{z, x\circ y\}
 =\{z ,x\}y'+(-1)^{|x||z|}x\{z,y'\}+z'(x\circ y) -z d(x\circ y) \\
 = (zx'-z'x)y' + (-1)^{|x||z|} xzy'' - z(xy')' + z'xy' \\
 =zx'y' - z'xy' + zxy'' - zxy'' - zx'y' +z'xy' = 0.
\end{multline*}
Hence, the entire differential ideal $I_V$ generated by $u(x,y)$ is invariant under the bracket 
defined by \eqref{eq:NovPoisson_s_Bracket}.
\end{proof}

\begin{theorem}\label{thm:Nov-Poisson}
Let $V$ be a Novikov superalgebra. 
Then the GD-algebra $V^{(-)}$ is special.
\end{theorem}

\begin{proof}
Lemmas \ref{lem:BrJacobi}, \ref{lem:Der_bracket},  \ref{lem:bracket_Ideal}
show the quotient $U(V) = s\Com\Der \<X_0\cup X_1, d\>/ I_V$
is a differential Poisson algebra. Theorem \ref{thm:NovToComDer}
guaranties that $V$ embeds into $U(V)$ with 
$u\circ v = ud(v)$, $u,v\in V$. Moreover, 
\eqref{eq:NovPoisson_s_Bracket} implies
$\{x,y\} = -x'y +xy' = x\circ y -(-1)^{|x||y|}y\circ x$
for $x,y\in X$. Therefore, $U(V)$ is a differential Poisson 
enveloping superalgebra of $V^{(-)}$.
\end{proof}

\section{Conformal superalgebras}\label{sec:Conformal}

In this section we recall the main definitions concerning Lie conformal (super)al\-geb\-ras 
and their relations to GD-(super)algebras.

Let $L$ be a Lie superalgebra over a field $\Bbbk $ of characteristic zero
and let $H=\Bbbk [\partial ]$ be the algebra of polynomials in one variable.

Two formal distributions 
$a(z),b(z)\in  L[[z,z^{-1}]]$ are said to be {\em mutually local}
if there exists $N\ge 0$ such that 
\[
 [a(w),b(z)](w-z)^N=0\in L[[z,z^{-1},w,w^{-1}]]. 
\]
For a pair of mutually local formal distributions $a(z)$, $b(z)$ 
their product $[a(w),b(z)]$ may be uniquely presented as a finite 
distribution 
\[
 [a(w),b(z)] = \sum\limits_{n=0}^{N-1} c_n(z) \dfrac{1}{n!} \dfrac{\partial^n \delta(w-z)}{\partial z^n},
\]
where $c_n(z)\in L[[z,z^{-1}]]$, $\delta(w-z) = \sum\limits_{s\in \mathbb Z} w^sz^{-s-1}$
is the formal delta-function. 

The collection of formal distributions $c_n(z)$, $n=0,1,\dots, N-1$, associated to a given pair $a(z)$, $b(z)$
may be written as a polynomial in a new formal variable $\lambda $:
\[
 [a(z)\oo\lambda b(z)] = \sum\limits_{n=0}^{N-1} \dfrac{\lambda^n}{n!} c_n(z) \in (L[[z,z^{-1}]])[\lambda ].
\]
The latter space may be identified with 
$\Bbbk [\partial,\lambda ]\otimes _H L[[z,z^{-1}]] $, where $\partial $ acts on $L[[z,z^{-1}]]$ 
as the ordinary derivation with respect to~$z$.

An $H$-invariant subspace $C$ of $L[[z,z^{-1}]]$ which consists of pairwize mutually local distributions 
and is closed with respect to the operation $[\cdot\oo\lambda \cdot ]$ (i.e., $a,b\in C$ implies $[a\oo\lambda b]\in C[\lambda ]$)
provides an example of the following class of algebraic structures.

\begin{definition}[{\cite[Chapter~2]{Kac1998}}]\label{defn:ConfAlgebra}
A $\mathbb Z_2$-graded $H$-module $C=C_0\oplus C_1$ equipped with a linear map ($\lambda$-{\em bracket})
\[
 [\cdot \oo{\lambda }\cdot ]: C\otimes C \to \Bbbk[\partial,\lambda]\otimes _H C\simeq C[\lambda ]
\]
is called a {\em Lie conformal superalgebra} if $[C_i\oo\lambda C_j]\subseteq C_{(i+j)\pmod 2} [\lambda ]$
and
\begin{gather}
{}
[\partial x\oo\lambda y ] = -\lambda [x\oo\lambda y], \label{eq:Sesqui1} \\ 
[x\oo\lambda \partial y ] = (\partial+\lambda )[x\oo\lambda y], \label{eq:Sesqui2} \\
[x\oo\lambda y] = (-1)^{|x||y|} [y\oo{-\partial-\lambda } x], \label{eq:AcommZ2} \\
[x\oo\lambda [y\oo\mu z]] - (-1)^{|x||y|}[y\oo\mu[x\oo\lambda z]] = [[x\oo\lambda y]\oo{\lambda+\mu} z]
\end{gather}
for all $x,y\in C$. Here $|x|$ stands for the parity of a homogeneous element $x\in C_0\cup C_1$.
\end{definition}

An operation $[\cdot\oo\lambda \cdot]$ satisfying \eqref{eq:Sesqui1} and \eqref{eq:Sesqui2}
is said to be {\em sesqui-linear}.

In particular, if $C$ is a free $H$-module 
with a basis $B$ then the 
$\lambda $-bracket on $C$ is uniquely determined by polynomials 
$p_{x,y}^z \in \Bbbk[\partial,\lambda ]$, 
$x,y,z\in B$, so that
\[
[x\oo\lambda y] = \sum\limits_{z\in B} p_{x,y}^z(\partial, \lambda )z.
\]

\begin{example}[\cite{Xu2000}]\label{exmp:Quadratic}
Let $V$ be a GD-superalgebra. Then the space of Laurent polynomials 
$V[t,t^{-1}]$ is a Lie superalgebra relative to the bracket
\[
 [at^n,bt^m] = [a,b]t^{n+m} +n (a\circ b) t^{n+m-1} - (-1)^{|a||b|} (b\circ a)t^{n+m-1}, 
 \quad a,b\in V.
\]
For $a\in V$, let $a(z)$ stand for the formal distribution
\[
 a(z) = \sum\limits_{n\in \mathbb Z} at^n z^{-n-1} \in (V[t,t^{-1}])[[z,z^{-1}]].
\]
The linear span of all such formal distributions 
and all their derivatives (with respect to $z$) is a Lie conformal superalgebra $L(V)$ isomorphic 
to $H\otimes V$ relative to the $\lambda $-bracket
\begin{equation}\label{eq:L-bracketGD}
  [a\oo\lambda b] = [a,b] + (-1)^{|a||b|}(\partial +\lambda)(b\circ a) + \lambda (a\circ b), 
  \quad a,b\in V.
\end{equation}
\end{example}

Conformal superalgebras obtained from GD-superalgebras as in Example 
\ref{exmp:Quadratic}
are called {\em quadratic Lie conformal superalgebras}.
As shown in \cite[Theorem~2.2]{Xu2000}, every Lie conformal superalgebra structure on 
a free $H$-module $L=H\otimes V$ given by linear polynomials 
$p_{x,y}^z(\partial, \lambda )$, $x,y,z\in V$, is quadratic.

An important particular example of a quadratic Lie conformal algebra comes 
from the 1-dimensional GD-algebra $V=\Bbbk v$ with $v\circ v = v$:
\[
[v\oo\lambda v] = \partial v + 2\lambda v.
\]
In this case, $L(V)$ is the {\em Virasoro conformal algebra}.

Suppose $V$ and $W$ are two $H$-modules.
The space of all linear maps 
 \[
  \alpha : V\to W[\lambda ], \quad v\mapsto \alpha_\lambda (v),
 \]
such that $\alpha_\lambda (\partial v) = (\partial+\lambda )\alpha_\lambda(v)$
for all
$v\in V$ is denoted by $\Chom (V,W)$ 
(the space of {\em conformal homomorphisms} from $V$ to~$W$) \cite{Kac1998}.
If $V=W$ then $\Chom(V,W)$ is denoted by $\Cend V$. 

If $V$ is a finitely generated $H$-module then the operation 
\[
 (\cdot \lambda \cdot): \Cend V\otimes \Cend V \to \Cend V[\lambda ]
\]
given by 
\begin{equation}\label{eq:L-prodCend}
 (\alpha \oo\lambda \beta)_\mu (v) = \alpha_\lambda (\beta_{\mu-\lambda }(v)), \quad v\in V,
\end{equation}
satisfies \eqref{eq:Sesqui1} and \eqref{eq:Sesqui2}.
Note that if $V$ is not a finitely generated $H$-module then we cannot say in general that 
$(\alpha\oo{\lambda }\beta )$ is a polynomial in~$\lambda $.

\begin{example}\label{gcN}
Let $U=U_0\oplus U_1$ be a finite-dimensional $\mathbb Z_2$-graded 
linear space, and let $V=H\otimes U = V_0\oplus V_1$
be the free $H$-module generated by~$U$, $V_i=H\otimes U_i$.
Then $\Cend V$ splits into the sum of even and odd components in 
a natural way:
\[
 (\Cend V)_0 = \Cend (V_0) \oplus \Cend (V_1), 
 \quad 
 (\Cend V)_1 = \Chom (V_0,V_1)\oplus \Chom(V_1,V_0),
\]
and the bracket 
\[
 [\alpha \oo\lambda \beta ] = (\alpha\oo{\lambda } \beta) - (-1)^{|\alpha||\beta|}
 (\beta \oo{-\partial -\lambda } \alpha) 
\]
turns $\Cend V$ into a Lie conformal superalgebra denoted $\gc V$.
\end{example}

If $V$ is a free $H$-module of rank $n+m$, 
where 
$n=\dim U_0$, $m=\dim U_1$, 
then $\gc V$ is denoted $\gc_{n|m}$. 
This conformal superalgebra may be presented as 
\[
 H\otimes H\otimes M_{n|m}(\Bbbk )\simeq H\otimes M_{n|m}(\Bbbk[x]),
\]
where $M_{n|m}$ stands for the $\mathbb Z_2$-graded algebra of $(n+m)$-matrices with 
even component 
$\begin{pmatrix} 
  M_n & 0 \\ 0 & M_m
 \end{pmatrix}$
 and odd component 
$\begin{pmatrix} 
  0 & M_{n,m} \\ M_{m,n} & 0
 \end{pmatrix}$.
The $\lambda $-bracket is given by 
\[
 [A(x)\oo\lambda B(x)] = A(x)B(x+\lambda ) - (-1)^{|A||B|} B(x)A(x-\partial-\lambda), 
\]
for homogeneous matrices $A,B\in M_{n|m}(\Bbbk[x])$.

A finite {\em representation} of a Lie conformal superalgebra $L$ 
is a homomorphism of Lie conformal superalgebras 
\[
 \rho : L\to \gc_{n|m}.
\]
In order to have a finite faithful representation (FFR), a Lie conformal superalgebra
has to be a torsion-free $H$-module since $\rho(\mathrm{tor_H} L)=0$. 

Every finite torsion-free Lie conformal superalgebra $L$ has a {\em regular} representation 
on itself: $\rho(a)=\alpha $ for $a\in L$, where $\alpha_\lambda (x) = [a\oo\lambda x]$.

In general, a representation of a Lie conformal superalgebra $L$ on a $\mathbb Z_2$-graded 
$H$-module $V=V_0\oplus V_1$ is a 
sesqui-linear map $\rho_\lambda (\cdot,\cdot): L\otimes V\to V[\lambda ]$
which respects the gradings and
\begin{equation}\label{eq:JacobiRepr}
\rho_\lambda (a,\rho_\mu(b,x)) - (-1)^{|a||b|} \rho_\mu (b,\rho_\lambda (a,x))  = \rho_{\lambda+\mu}([a\oo\lambda b],x),
\end{equation}
for $a,b\in L_0\cup L_1$, $x\in V$.

This is unknown whether all torsion-free finite Lie conformal (super)algebras 
have FFR. In the non-graded case, it is known \cite{Kol2016} that if 
the solvable radical of $L$ splits  then $L$ has a FFR. The proof of the latter 
result essentially involves representation theory of Lie conformal algebras 
\cite{ChengKac}. 
For example, if $V$ is the exceptional GD-algebra from Example~\ref{exmp:Heisenberg3}
then $L(V)$ has a split solvable radical. 
However, a quadratic Lie conformal (super)algebra may have a non-split 
solvable radical. For example, the Virasoro Lie conformal algebra has 
a non-split extension \cite[Theorem 7.2]{BKV99} 
corresponding to the 2-dimensional Novikov algebra $V=\Bbbk v+\Bbbk u$, 
where 
\[
 v\circ v = v+u,\quad v\circ u=0,\quad u\circ v = u,\quad u\circ u=0.
\]
Considered as an Abelian Lie algebra, $V$ is a GD-algebra that 
gives rise to a quadratic Lie conformal algebra $L(V)$ with 
a non-split solvable radical $H\otimes \Bbbk u$. 

In the next section, we will show that for every special GD-superalgebra 
$V$ the corresponding Lie conformal superalgebra $L(V)$ has a FFR.

\section{Poisson conformal superalgebras}\label{sec:Poisson}

In the study of Ado-type problems for Lie conformal algebras,
Poisson structures play an important role.

\begin{definition}\label{defn:PoissonConf}
Let $P$ be a $\mathbb Z_2$-graded $H$-module endowed with two sesqui-linear operations 
\[
\begin{gathered}
 {}
 [\cdot \oo\lambda \cdot ]: P\otimes P\to P[\lambda ], \\
 (\cdot \oo\lambda \cdot ): P\otimes P\to P[\lambda ],
\end{gathered}
\]
which respect the grading and satisfy the following conditions:
\begin{enumerate}
 \item $P$ is a Lie conformal superalgebra relative to $[\cdot \oo\lambda \cdot]$; 
 \item\label{cnd:Confass} $(a\oo\lambda (b\oo\mu c)) = ((a\oo\lambda b)\oo{\lambda +\mu } c)$
 for $a,b,c\in P$; 
 \item $(a\oo\lambda b) = (-1)^{|a||b|} (b\oo{-\partial-\lambda} a)$ 
for homogeneous $a,b\in P$.
 \end{enumerate}
Then $P$ is said to be a {\em Poisson conformal superalgebra}
if the following conformal analogue of the Leibniz rule holds:
\begin{equation}\label{eq:ConformalLeibniz}
 {} [a\oo\lambda (b\oo\mu c)] = ([a\oo\lambda b]\oo{\lambda+\mu} c) + (-1)^{|a||b|}(b\oo\mu [a\oo\lambda c]),
\end{equation}
for $a,b\in P_0\cup P_1$, $c\in P$.
\end{definition}

A simplest example of a Poisson conformal (super)algebra is provided by the current functor. Namely, 
if $\mathfrak p$ is an ordinary commutative (super)algebra with a (super)Poisson bracket $\{\cdot,\cdot\}$ 
then $\Cur \mathfrak p = H\otimes \mathfrak p\simeq \mathfrak p[\partial ]$ 
equipped with 
\[
\begin{gathered}
 (a(\partial )\oo\lambda b(\partial )) = a(-\lambda )b(\partial+\lambda ), \\
 [a(\partial )\oo\lambda b(\partial )] = \{a(-\lambda ),b(\partial+\lambda )\}
\end{gathered}
\]
is a Poisson conformal (super)algebra.

An associative conformal algebra \cite{Kac1998} 
is an $H$-module $A$ equipped with a sesqui-linear operation $(\cdot\oo\lambda \cdot)$
satisfying the Condition \ref{cnd:Confass} of Definition \ref{defn:PoissonConf}. 
Assuming $A$ is a $\mathbb Z_2$-graded associative conformal algebra, 
the new (commutator) operation 
\[
 [a\oo\lambda b]=(a\oo\lambda b) - (-1)^{|a||b|} (b\oo{-\partial-\lambda } a),
 \quad 
a,b\in A_0\cup A_1,
\]
turns the $H$-module $A$ into a Lie conformal superalgebra $A^{(-)}$ (see \cite[p.~323]{Ro2000} for non-graded case). 

Given a Lie conformal superalgebra $L$ and a $\mathbb Z_2$-graded associative conformal algebra $A$, 
we say $A$ is an {\em associative conformal envelope} of $L$ if there exists a 
homogeneous homomorphism $\tau : L\to A^{(-)}$ of conformal algebras such that $A$ is generated 
(as an associative conformal algebra) by the image of~$L$.
For a fixed $L$, there exists a lattice of universal associative conformal envelopes of $L$ 
corresponding to different associative locality bounds on the elements of $L$ (see \cite[Section~6]{Ro2000}).
It may happen that neither of these universal envelopes contains the isomorphic image of $L$, i.e., 
there exist Lie conformal (super)algebras that cannot be embedded into associative ones.

Suppose $A$ is an associative conformal envelope of a Lie conformal superalgebra $L$. 
Then $A$ has a natural filtration as an $H$-module:
\[
 0=A_0\subset A_1\subset A_2\subset \dots , 
\]
where $ A_1=\tau(L)$, $A_{n+1} = A_n+H\{ (\tau(L)\oo{\lambda } A_{n})|_{\lambda =\alpha} : \alpha \in \Bbbk \}$.
Then 
\[
 (A_n\oo\lambda A_m) \subseteq A_{n+m}[\lambda ], \quad [A_n\oo\lambda A_m] \subseteq A_{n+m-1}[\lambda ], 
\]
so the associated graded $H$-module $\mathrm{gr}\,A$  has a well-defined structure 
of a Poisson conformal superalgebra.

For example, {\em conformal Weyl algebra} $\Cend_{1,x}$ (see \cite{BKL2003}) is an associative 
envelope of the Virasoro conformal algebra. The corresponding Poisson conformal algebra 
$\mathrm{gr}\,\Cend_{1,x} \simeq \Bbbk [\partial ]\otimes x\Bbbk[x]$ has the following 
operations:
\[
 (x^n\oo\lambda x^m) = x^{n+m},\quad [x^n\oo\lambda x^m] = (n\partial + (n+m)\lambda )x^{n+m-1}.
\]
More examples of Poisson conformal superalgebra structures on universal associative envelopes of 
Lie conformal superalgebras can be found in \cite{Kol:PoisPrepr}.

Hereinafter, we will need the following example of a Poisson conformal superalgebra. 

\begin{lemma}\label{lem:QuadCurrent}
Let $\mathfrak p = \mathfrak p_0 \oplus \mathfrak p_1$ be an ordinary Poisson superalgebra equipped with an even derivation $d: a\mapsto a'$. 
The latter means $d(\mathfrak p_i)\subseteq \mathfrak p_i$, $i=0,1$,
$(ab)'=ab'+a'b$, $\{a,b\}' = \{a',b\} + \{a,b'\}$ for $a,b\in \mathfrak p$.
Then $L(\mathfrak p, d) = H\otimes \mathfrak p$ equipped with operations
\[
 (a\oo{\lambda } b) = ab,\quad 
 [a\oo\lambda b] = \{a,b\} + \partial a'b + \lambda (ab)', 
\]
for $a,b\in \mathfrak p_0\cup \mathfrak p_1$, 
is a Poisson conformal superalgebra.
\end{lemma}

\begin{proof}
The operation $(\cdot \oo\lambda \cdot)$ is obviously associative and (super-) commutative. 
By definition,
$[\cdot \oo\lambda \cdot ]$ is exactly the quadratic Lie conformal bracket 
on the GD-superalgebra obtained from $\mathfrak p$ relative to~$d$. 
The Leibniz rule 
\eqref{eq:ConformalLeibniz} is straightforward to check. On the one hand,
\begin{multline}\nonumber
 [a\oo\lambda (b\oo\mu c)] - (-1)^{|a||b|} (b\oo\mu [a\oo\lambda c]) \\
 =
 \{a,bc\} + \partial a'bc +\lambda (abc)' 
 - (-1)^{|a||b|} \big (
 b\{a,c\} + (\partial+\mu)ba'c  +\lambda b(ac)' 
     \big ) \\
 =
 \{a,b\}c - \mu a'bc +\lambda ab'c
\end{multline}
for homogeneous $a,b,c\in \mathfrak p$.
On the other hand, 
\[
 ([a\oo\lambda b]\oo{\lambda +\mu} c)
 =
 \{a,b\}c - (\lambda +\mu)a'bc + \lambda (ab)'c
 =
 \{a,b\}c - \mu a'bc +\lambda ab'c.
\]
Hence, $L(\mathfrak p,d)$ is a Poisson conformal superalgebra.
\end{proof}

Suppose $L$ is a Lie conformal superalgebra with a representation $\rho $ on an $H$-module $V$. 
Then a {\em deformation} of $\rho $ is a representation $\rho_\varepsilon$ 
of $L$ on the $H$-module $V\oplus V\varepsilon\simeq \Bbbk[\varepsilon]\otimes V/(\varepsilon^2) $, 
where 
\[
 \rho_\varepsilon (a) = \alpha^{(\varepsilon )}, 
 \quad 
 \alpha^{(\varepsilon)}_\lambda (x) = \alpha_\lambda (x)  +\varepsilon \varphi_\lambda (a,x),
\]
for $a\in L$, $x\in V$, $\alpha = \rho(a)$.
Here $\varphi_\lambda(\cdot,\cdot) : L\otimes V \to V[\lambda ]$ 
is a sesqui-linear map which has to satisfy the following equation (a consequence of the 
Jacobi identity \eqref{eq:JacobiRepr}):
\begin{multline}\label{eq:CendCocycle}
\varphi_{\lambda +\mu} ([a\oo\lambda b],x) 
=
\varphi_\lambda (a,\rho_\mu(b,x))  + \rho_\lambda (a, \varphi_\mu(b,x)) \\
-(-1)^{|a||b|} \varphi_\mu (b, \rho_\lambda (a,x))-(-1)^{|a||b|} \rho_\mu (b,\varphi_\lambda (a,x)),
\end{multline}
for $a,b\in L_0\cup L_1$, $x\in V$.
In the case when $V$ is a finitely generated $H$-module, the relation 
\eqref{eq:CendCocycle}
exactly means that $\varphi$ is a 1-cocycle in $Z^1(L,\Cend V)$ (see \cite{BKV99}).

\begin{proposition}\label{prop:PoissCocycle}
Let $L$ be a graded Lie conformal subalgebra in a Poisson conformal superalgebra $P$. 
Then $L$ has a regular representation on $P$, and 
the sesqui-linear map $\varphi_\lambda (\cdot,\cdot ): L\otimes P \to P[\lambda ]$
given by 
\[
 \varphi_\lambda (a,x) = \lambda (a\oo\lambda x), \quad a\in L,\ x\in P,
\]
satisfies \eqref{eq:CendCocycle}.
\end{proposition}

\begin{proof}
Obviously, the map $\rho_\lambda (a,x) = [a\oo\lambda x] $, $a\in L$, $x\in P$, 
is a representation of $L$ on $P$. 
Then the right-hand side of \eqref{eq:CendCocycle} 
can be transformed by \eqref{eq:ConformalLeibniz} as
\begin{multline}
\lambda (a\oo\lambda [b\oo\mu x])  + \mu[a\oo\lambda (b\oo\mu x)] \\
-(-1)^{|a||b|} \mu (b\oo\mu [a\oo\lambda x])-(-1)^{|a||b|} \lambda [b\oo\mu (a\oo\lambda x)] \\
=
\lambda ([a\oo\lambda b] \oo{\lambda +\mu } x) + \mu ([a\oo\lambda b] \oo{\lambda +\mu } x) 
=\varphi_{\lambda +\mu }([a\oo\lambda b],x).
\end{multline}
\end{proof}

\begin{corollary}[\cite{Kol:PoisPrepr}]\label{cor:TwistRepr}
 Let $L$ be a graded Lie conformal subalgebra of a Poisson conformal superalgebra. 
 Then the map 
 \[
  \hat \rho_\lambda (a,x) = [a\oo\lambda x] + \lambda (a\oo\lambda x), \quad a\in L,\ x\in P,
 \]
is a representation of $L$ on $P$.
\end{corollary}

\begin{proof}
 In order to check the Jacobi identity \eqref{eq:JacobiRepr} for $\hat \rho$, note that 
 the desired equation
 \[
  \hat\rho_\lambda (a,\hat\rho_\mu(b,x)) - (-1)^{|a||b|} \hat\rho_\mu (b,\hat\rho_\lambda (a,x))  = \hat\rho_{\lambda+\mu}([a\oo\lambda b],x)
 \]
splits into three equations: the first one is exactly the Jacobi identity for regular representation $\rho_\lambda (a,x)=[a\oo\lambda x]$, 
the second one is \eqref{eq:CendCocycle} for $\varphi_\lambda (a,x)=\lambda (a\oo\lambda x)$, 
and the third one is 
\[
 \lambda \mu (a\oo\lambda (b\oo\mu x)) - (-1)^{|a||b|}\lambda \mu (b\oo\mu (a\oo\lambda x))=0
\]
which also holds due to conformal commutativity and associativity of $(\cdot \oo\lambda \cdot )$:
\begin{multline}\nonumber
 (a\oo\lambda (b\oo\mu x)) 
 = ((a\oo\lambda b)\oo{\mu+\lambda } x)
 = (-1)^{|a||b|} ((b\oo{-\partial -\lambda } a)\oo{\lambda +\mu } x) \\
 =(-1)^{|a||b|} ((b\oo{\mu } a)\oo{\lambda +\mu} x)
 =(-1)^{|a||b|} (b\oo{\mu } (a\oo{\lambda } x)).
\end{multline}
\end{proof}

\begin{corollary}
 If for a finite Lie conformal superalgebra $L$
 there exists an embedding $\tau:L\to P$ of $L$ into a Poisson conformal 
 superalgebra $P$ in such a way that 
 $(\tau(a)\oo\lambda \tau(L))\ne 0$ for all $0\ne a\in L$
 then 
 $L$ has a FFR.
\end{corollary}

\begin{proof}
 Indeed, if we consider $V=L$ as a regular $L$-module, 
 $M = \tau(L)+\Bbbk[\partial] \otimes \mathrm{span}\,\{(\tau(L)\oo\lambda \tau(L))|_{\lambda =\alpha }: \alpha\in \Bbbk \}$
 as an $L$-submodule of $P$, and 
 \[
  \langle \cdot\oo\lambda \cdot\rangle : L\otimes V \to M[\lambda ]
 \]
given by $\langle a\oo\lambda b\rangle = (\tau (a)\oo\lambda \tau(b))$, 
then all conditions of \cite[Theorem~3]{Kol11} hold and $L$ has a FFR.
\end{proof}

In particular, if $L$ satisfies the Poincar\'e--Birkhoff--Witt condition \cite{Ro2000} then one may choose $P$ to be the 
associated graded conformal algebra of the appropriate universal associative conformal envelope of $L$. 
Therefore, for conformal algebras the PBW theorem implies 
the Ado Theorem immediately \cite{Kol11}.

\begin{theorem}\label{thm:QuadSpecFFR}
 Let $V$ be a finite-dimensional special GD-superalgebra. Then 
 the Lie conformal superalgebra $L(V)$ has a finite faithful 
 conformal representation.
\end{theorem}

\begin{proof}
Let us fix linear bases $X_0$ and $X_1$ of $V_0$ and $V_1$, respectively, and let $X=X_0\cup X_1$.
A special GD-superalgebra embeds into its universal enveloping differential Poisson superalgebra 
which can be constructed as follows. 
Denote by $F = s\Pois\Der\<X_0\cup X_1, d\>$ the free differential Poisson superalgebra 
with an identity element $1$ (generated by  
even elements $X_0$, odd elements $X_1$) and an even derivation~$d$. 
Let $I_V$ stand for the (differential) ideal of $F$
generated by 
\begin{equation}\label{eq:IV_generator}
xd(y) - x\circ y, \quad x,y\in X,
\end{equation}
 where $x\circ y$ is a linear form in $\Bbbk X$ representing the Novikov product in~$V$.
The quotient $F/I_V$ is the universal 
 enveloping differential Poisson superalgebra for $V$ denoted $P_d(V)$. 
If $V$ is a special GD-superalgebra then $V$ embeds into $P_d(V)$.

The free algebra $F$ may be presented as 
\[
F = \bigoplus_{n\in \mathbb Z} F_n,
\]
where $F_n$ consists of all elements of weight $n$. Recall that the weight in a free differential 
Poisson (super)algebra 
generated by a set $X$ is defined as follows \cite{KSO2019}:
\[
\begin{gathered}
 \wt x = -1\text{ for } x\in X, \quad \wt(1)=0,\\
 \wt (uv) = \wt u + \wt v, \quad \wt \{u,v\} = \wt u + \wt v +1, \\
 \wt d(u)= \wt u +1.
\end{gathered}
\]
Since all elements in \eqref{eq:IV_generator} are $\wt$-homogeneous, the ideal $I_V$ is $\wt$-homogeneous 
and the algebra $P_d(V)$ inherits the grading:
\[
P_d(V) = \bigoplus_{n\in \mathbb Z} U_n, \quad U_n = F_n/I_V\cap F_n.
\]
Note that $V\simeq U_{-1}$. The latter was shown in \cite[Theorem~10]{KSO2019} for non-graded case, 
the same reasonings work for superalgebras.
 
Lemma \ref{lem:QuadCurrent} states that $L(P_d(V), d)=H\otimes P_d(V)$ is a Poisson conformal superalgebra. 
Then by Corollary \ref{cor:TwistRepr} the Lie conformal superalgebra
 $L=L(V)$ has a representation on $M = L(P_d(V), d)$ given by
\begin{equation}\label{eq:ReprPoisTwist}
\hat \rho_\lambda (a,u) = \{a,u\} +\partial d(a)u + \lambda (d(au)+au), \quad a\in V, \ u\in P_d(V).
\end{equation}
Obviously, for every $m\in \mathbb Z$ the space 
$M_{\le m} = H\otimes \bigoplus_{n\le m} U_n$
is  a conformal $L$-submodule of~$M$. 
In particular, $\bar M = M_{\le 0}/M_{\le -2} \simeq H\otimes (U_{-1}\oplus U_0)$
is a conformal $L$-module corresponding to a representation $\bar\rho $
defined by \eqref{eq:ReprPoisTwist} for $u\in U_0$ and $\bar\rho_\lambda (a,u)= [a\oo\lambda u]$ 
for $u\in U_{-1}\simeq V$.
It is easy to see that the representation of $L$ on $\bar M$ is faithful:
$\bar \rho_\lambda (a,1) = (\partial+\lambda )d(a) + \lambda a \ne 0$ for $a\in V$, $a\ne 0$.
However, it is not yet finite in general since $\dim U_0$ may not be finite.

Note that for every $a\in V\simeq U_{-1}$ the map $\mu_a: u\mapsto au$ maps $U_0$ to $U_{-1}$. 
Since $\dim U_{-1}=\dim V<\infty $, the intersection of all $\Ker\mu_a$, $a\in V$, is a subspace 
of finite codimension. So, consider $N = \{u\in U_0 \mid Vu=0 \}$. For every $u\in N$ we have 
\[
\bar\rho_\lambda (a,u) = \{a,u\} + \partial d(a)u, \quad a\in V.
\]
Given $b\in V$, 
$b\{a,u\} = (-1)^{|a||b|}(\{a,bu\} + \{a,b\}u) = 0$, 
$b d(a)u = (-1)^{|a||b|} d(a) bu = 0$.
Therefore, $H\otimes N$ is a conformal $L$-submodule of $\bar M$.
Finally,
\[
\bar M/N \simeq H\otimes (U_{-1}\oplus U_0/N)
\]
is a finite faithful conformal $L$-module.
\end{proof}

\begin{corollary}
If a GD-(super)algebra $V$ is constructed from a Novikov (super)algebra  
with respect to the commutator then the 
corresponding quadratic Lie conformal superalgebra $L(V)$ has a FFR.
\end{corollary}


\begin{thebibliography}{99}

\bibitem{BaiDeng01}
C.~Bai, D.~Meng, 
The classification of Novikov algebras in low dimensions,
J. Phys. A, Math. Gen. {\bf 34} (2001) (8) 1581--1594.


\bibitem{BDK}
B. Bakalov, A. D’Andrea, V. G. Kac,
Theory of finite pseudoalgebras, 
Adv. Math. {\bf 162} (2001) 1--140.

\bibitem{BKV99}
B. Bakalov, V. G. Kac, A. Voronov,
Cohomology of conformal algebras,
Comm. Math. Phys.  {\bf 200} (1999) 561--589.

\bibitem{BN85}
A. A. Balinskii, S. P. Novikov, 
Poisson brackets of hydrodynamic type, Frobenius algebras and Lie algebras, 
Sov. Math. Dokl. {\bf 32} (1985) 228--231.

\bibitem{BCZ-2017}
L.A. Bokut, Y. Chen, Z. Zhang, Gr\"obner--Shirshov bases method for Gelfand--Dorfman--Novikov algebras, 
J. Algebra Appl. 16 (1) (2017) 22.

\bibitem{BKL2003}
C. Boyallian, V. G. Kac, J.-I. Liberati, 
On the classification of subalgebras of $\Cend_N$ and $\mathrm{gc}_N$.
J.~Algebra {\bf 260} (2003) 32--63.

\bibitem{BdG2013}
D. Burde, W. de Graaf, 
Classification of Novikov algebras,
Appl. Algebra Eng. Commun. Comput. {\bf 24} (2013) (1) 1--15.

\bibitem{ChengKac}
S.-J. Cheng, V. G. Kac, 
Conformal modules,
 Asian J. Math.  {\bf 1} (1997) 181--193.

\bibitem{DzhLowf} 
A. S. Dzhumadil’daev, C. L\"ofwall, 
Trees, free right-symmetric algebras, free Novikov algebras and identities, 
Homology, Homotopy Appl. {\bf 4} (2002) (2) 165--190.

\bibitem{GuoKeig08}
Li Guo, W. Keigher,
On differential Rota-Baxter algebras,
J. Pure Appl. Algebra {\bf 212} (2008) (3) 522--540. 

\bibitem{GD83} 
I.~M.~Gelfand, I.~Ya.~Dorfman, 
Hamilton operators and associated algebraic structures, 
Functional analysis and its application {\bf 13} (1979) (4) 13--30.

\bibitem{GK94}
V. Ginzburg, M. Kapranov,
Koszul duality for operads,
Duke Math. J. {\bf 76} (1994) (1) 203--272.

\bibitem{HongWu17}
Y. Hong, Z. Wu,
Simplicity of quadratic Lie conformal algebras,
Comm. Algebra {\bf 45} (2017) (1) 141--150. 

\bibitem{Kac1998}
V. G. Kac, 
Vertex Algebras for Beginners. University Lecture Series, {\bf 10}, 2nd edn. 
American Mathematical Society, Providence, 1996 (1998).

\bibitem{Kol11}
P. S. Kolesnikov, 
On finite representations of conformal algebras.
J. Algebra {\bf 331} (2011) 169--193.

\bibitem{Kol2016}
P. S. Kolesnikov,
The Ado theorem for finite Lie conformal algebras with Levi decomposition.
J. Algebra Appl. {\bf 15} (2016) Art.no.~1650130.

\bibitem{KSO2019}
P. S. Kolesnikov, B. Sartayev, A. Orazgaliev,
Gelfand--Dorfman algebras, derived identities, and the Manin product of operads,
Journal of Algebra {\bf 539} (2019) 260--284.


\bibitem{Kol:PoisPrepr}
P. S. Kolesnikov, 
Universal enveloping Poisson conformal algebras,
Internat. J. Algebra Comput.  {\bf 30} (2020) no.~5, 1015--1034.

\bibitem{Lod10}
J.-L. Loday,
On the operad of associative algebras with derivation,
Georgian Math. J. {\bf 17} (2010) (2) 347--372. 

\bibitem{LodVal08}
J.-L.~Loday, B.~Vallette,
 Algebraic Operads, 
Gr\"undlehren der mathematischen Wissenschaften {\bf 346}, 
Springer-Verlag, 2012.

\bibitem{MikhSh10}
A. A. Mikhalev, I. P. Shestakov, 
PBW-pairs of varieties of linear algebras,
Comm. Algebra {\bf 42} (2014) 667--687.

\bibitem{Osborn}
J. M. Osborn, 
Novikov algebras, 
Nova J. Algebra \& Geom. {\bf 1} (1992) 1--14.

\bibitem{Ro2000}
M. Roitman,
Universal enveloping conformal algebras, 
{Sel. Math., New Ser.} {\bf 6} (2000) 319--345.

\bibitem{Shaf}
I. R. Shafarevich, 
Degeneration of semisimple algebras,
Comm. Algebra {\bf 29} (2001) no.~9, 3943--3960. 

\bibitem{Sh58}
A. I. Shirshov, 
On free Lie rings, 
Mat. Sb., 45(87) (1958) 113--122.

\bibitem{Xu1996}
X. Xu, 
On simple Novikov algebras and their irreducible modules,
J. Algebra {\bf 185} (1996) (3) 905--934.

\bibitem{Xu2000} 
X. Xu, 
Quadratic Conformal Superalgebras, 
J. Algebra {\bf 231} (2000) 1--38.

\bibitem{Xu2001}
X. Xu, 
Classification of simple Novikov algebras and their irreducible modules of characteristic~0,
J. Algebra {\bf 246}, (2001) (2) 673--707.

\bibitem{Xu2002}
X. Xu, 
Gel'fand--Dorfman bialgebras,
Southeast Asian Bull. Math. {\bf 27} (2003) 561--574.


\bibitem{Zelm}
E. I. Zel'manov,
On a class of local translation invariant Lie algebras,
Sov. Math. Dokl. {\bf 35} (1987) 216--218.

\bibitem{ZCB2019} Z. Zhang, Y. Chen, L.A. Bokut,
Free Gelfand--Dorfman--Novikov superalgebras and a Poincar\'e--Birkhoff--Witt type theorem,
Int. J. Algebra Comput. 29 (2019), no. 3, 481--505.

\end{thebibliography}
\end{document}